\documentclass[12pt]{amsart}

\usepackage[
    colorlinks=true,       
    linkcolor=magenta,          
    citecolor=magenta,        
    filecolor=magenta,      
    urlcolor=magenta,           
    ]{hyperref}

\usepackage{xcolor}
\usepackage{amsmath,amsthm,amscd,amssymb,amsfonts, amsbsy}
\usepackage{latexsym}
\usepackage{exscale}
\usepackage{mathrsfs,amssymb}
\usepackage{enumerate}

\newtheorem{theorem}[equation]{Theorem}

\newtheorem{cor}[equation]{Corollary}
\newtheorem{lemma}[equation]{Lemma}

\newtheorem*{Remark}{Remark}
\newtheorem{definition}[equation]{Definition}

\numberwithin{equation}{section}

\def\0{{\rm \bf{0}}}

\begin{document}

\title{Weighted mixed weak-type inequalities for multilinear fractional operators}

\author{M. Bel\'en Picardi}
\address{Departamento de Matem\'atica\\
Universidad Nacional del Sur\\
Bah\'ia Blanca, 8000, Argentina}\email{belen.picardi@uns.edu.ar}

\begin{abstract}
The aim of this paper is to obtain mixed weak-type inequalities for multilinear fractional operators, extending results by F. Berra, M. Carena and G. Pradolini \cite{BCP}. We prove that, under certain conditions on the weights, there exists a constant $C$ such that
$$\Bigg\| \frac{\mathcal G_{\alpha}(\vec f \,)}{v}\Bigg\|_{L^{q, \infty}(\nu v^q)} \leq C \ \prod_{i=1}^m{\|f_i\|_{L^1(u_i)}},$$ where $\mathcal G_{\alpha}(\vec f \,)$ is the multilinear maximal function $\mathcal M_{\alpha}(\vec f\,)$ that was introduced by K. Moen in \cite{M} or the multilineal fractional integral $\mathcal I_{\alpha}(\vec f \,)$. As an application a vector-valued weighted mixed inequality for $\mathcal I_{\alpha}(\vec f \,)$ will be provided as well.
\end{abstract}


\maketitle

\section{Introduction}
E. Sawyer \cite{S} proved in 1985 the following mixed weak-type inequality.

\begin{theorem}[\cite{S}]\label{sawyer}
If $u,v \in A_1$, then there is a constant $C$ such that for all $t>0$,
\begin{equation}\label{saw}
uv\Big\{ x \in \mathbb{R}  :\; \frac{M(fv)(x)}{v(x)} \; > t \Big\}
 \; \leq \; \frac{C}{t} \int_{\mathbb{R}} |f(x)|u(x)v(x)\,dx.
\end{equation}
\end{theorem}

This estimate is a highly non-trivial extension of the classical weak type $(1,1)$ inequality for the maximal operator due to the presence of the weight function $v$ inside the distribution set. Note that if $v=1$, this result is a well known estimate due to C. Fefferman - E. Stein \cite{FS}. (\ref{saw}) also holds, if $u\in A_1$ when $v\in A_1$, see \cite{LeOP}.

In 2005, D. Cruz-Uribe, J. M. Martell and C. P\'erez \cite{CMP} extended (\ref{saw}) to $\mathbb{R}^n$. Furthermore, they settled that estimate for Calder\'on-Zygmund operators, answering affirmatively and extending a conjecture raised by E. Sawyer for the Hilbert transform \cite{S}. The precise statement of their result is the following.

\begin{theorem}[\cite{CMP}]\label{cmpteo}
If $u,v \in A_1$, or $u\in A_{1}$ and $uv\in A_{\infty}$, then there is a constant $C$ such that for all $t>0$,
\begin{equation}\label{sawcmp}
uv\Big\{ x \in \mathbb{R}^n  :\; \frac{|T(fv)(x)|}{v(x)} \; > t \Big\}
 \; \leq \; \frac{C}{t} \int_{\mathbb{R}^n} |f(x)|u(x)v(x)\,dx,
\end{equation} 
where $T$ is a Calder\'on-Zygmund operator with some regularity.
\end{theorem}

Quantitative versions of the previous result were obtained in \cite{OPR} and also some counterparts for commutators in \cite{BCP1}.

In \cite{CMP}, D. Cruz-Uribe, J.M. Martell and C. P\'erez conjectured that (\ref{sawcmp}) and (\ref{saw}) should hold for $v\in A_{\infty}$. This result is the most singular case, due to the fact that the $A_{\infty}$ condition is the weakest possible asumption within the $A_p$ classes.

Recently K. Li, S. Ombrosi and C. P\'erez \cite{LOP1} solved that conjecture. They proved the following theorem.

\begin{theorem}[\cite{LOP1}]
Let $v\in A_{\infty}$ and $u\in A_1$. Then there is a constant $C$ depending on the $A_1$ constant of $u$ and the $A_{\infty}$ constant of $v$ such that
$$\Bigg\| \frac{T(fv)}{v}\Bigg\|_{L^{1, \infty}(uv)} \leq C \ \|f\|_{L^{1,\infty}(uv)} $$ where $T$ can be the Hardy-Littlewood maximal function, any Calder\'on-Zygmund operator or any rough singular integral.
\end{theorem}

\bigskip

In 2009, Lerner \textit{et al.} \cite{LOPTT} introduced the multi(sub)linear maximal function $\mathcal{M}$ defined by $$\mathcal M(\vec f\,)(x)=\sup_{Q\ni x
}\prod_{i=1}^m\frac{1}{|Q|}\int_Q|f_i(y_i)|dy_i,$$ where $\vec{f}=(f_1,...,f_m)$ and the supremum is taken over all cubes $Q$ containing $x$.

This maximal operator is smaller than the product $\prod_{i=1}^m Mf_i$, which was the auxiliar operator used previously to estimate multilinear singular integral operators.

\bigskip

There is a conection between multilinear operators and mixed weak-type inequalities (see \cite{LOPTT} or \cite{LOP1}). In fact, in a recent joint work with K. Li and S. Ombrosi \cite{LOP} we proved the following theorem.

\begin{theorem}[\cite{LOP}]\label{muczo}
Let $\mathcal{T}$ be a multilinear Calder\'on-Zygmund operator, $\vec{w}=(w_1,...,w_m)$ and $\nu = w_1^\frac{1}{m}...w_m^\frac{1}{m}$. Suppose that  $\vec{w} \in A_{(1,...1)}$ and $\nu v^{\frac 1m}\in A_\infty$ or  $w_1,...,w_m \in A_1$ and $v \in A_\infty$.  Then there is a constant $C$ such that
$$\Bigg\| \frac{\mathcal{T}(\vec f\,)(x)}{v}\Bigg\|_{L^{\frac{1}{m}, \infty}(\nu v^\frac{1}{m})} \leq C \ \prod_{i=1}^m{\|f_i\|_{L^1(w_i)}}.$$
\end{theorem}

\bigskip
We remit the reader to Section 2 for the definition of $A_{(1,...,1)}$ and more details about $A_{\vec{p}}$ weights.
\bigskip

The study of fractional integrals and associated maximal functions is important in harmonic analysis. We recall that the fractional integral operator or Riesz potential is defined by
$$I_{\alpha}f(x)=\int_{\mathcal{R}^n}{\frac{f(y)}{|x-y|^{n-\alpha}}\, dy}, \, \, \, 0<\alpha<n,$$ and the fractional maximal function by
$$M_{\alpha}f(x)=\sup_{x\in Q}{\frac{1}{|Q|^{1-\frac{\alpha}{n}}} \int_Q {|f(y)|dy}}, \, \, \, 0\leq \alpha <n,$$ where the supremum is taken over all cubes $Q$ containing $x$. Note that in the case $\alpha=0$ we recover the Hardy-Littlewood maximal operator. Properties of these operators can be found in the books by Stein \cite{Stein} and Grafakos \cite{G}.

F. Berra, M. Carena and G. Pradolini \cite{BCP} proved the following mixed weak-type inequality.
\begin{theorem}[\cite{BCP}]\label{teo1.7}
Let $0<\alpha <n$, $1\leq p <\frac{n}{\alpha}$ and $q$ satisfying $\frac{1}{q}=\frac{1}{p}-\frac{\alpha}{n}$. If $u$, $v$ are weights such that $u, v^{\frac{q}{p}} \in A_1$ or $uv^{\frac{-q}{p'}} \in A_1$ and $v\in A_{\infty}(uv^{\frac{-q}{p'}})$, then there exists a positive constant $C$ such that for every $t>0$
$$uv^{\frac{q}{p}}\Big\{ x \in \mathbb{R}^n  :\; \frac{I_{\alpha}(fv)(x)}{v(x)} > t \Big\}^{\frac{1}{q}}\leq \frac{C}{t} \Big( \int_{\mathbb{R}^n} {|f(x)|^p u(x)^{\frac{p}{q}} v(x) dx} \Big)^{\frac{1}{p}}$$ where $I_{\alpha}$ is the fractional integral or the fractional maximal function.
\end{theorem}

\bigskip

In the multilinear setting a natural way to extend fractional integrals is the following.

\begin{definition}
Let $\alpha$ be a number such that $0<\alpha<mn$ and $\vec{f}=(f_1,...,f_m)$ be a collection of functions on $\mathbb{R}^n$. We define the multilinear fractional integral as
$$\mathcal{I}_{\alpha}{\vec{f}(x)}=\int_{(\mathbb{R}^n)^m}{\frac{f_1(y_1)...f_m(y_m) \ d\vec{y}}{(|x-y_1|+...+|x-y_m|)^{mn-\alpha}}}$$
\end{definition}

K. Moen \cite{M} introduced the multi(sub)linear maximal operator $\mathcal{M}_{\alpha}$ asociated to the multilinear fractional integral $\mathcal{I}_{\alpha}$.

\begin{definition}
For $0\leq \alpha < mn$ and $\vec{f}=(f_1,...,f_m)$, we define the multi(sub)linear maximal operator $\mathcal{M}_{\alpha}$ by $$\mathcal{M}_{\alpha}\vec{f}(x)=\sup_{x\in Q}{\prod_{i=1}^m{\Bigg(\frac{1}{|Q|^{1-\frac{\alpha}{nm}}}\int_Q{|f_i(y_i)|dy_i}\Bigg)}}$$
\end{definition}

Observe that the case $\alpha=0$ corresponds to the multi(sub)linear maximal function $\mathcal{M}$ studied in \cite{LOPTT}.

\bigskip

At this point we present our contribution. Our first result is a counterpart of Theorem \ref{teo1.7} for multilinear fractional maximal operators.

\bigskip
\begin{theorem}\label{Max}
Let $0\leq \alpha<mn$. Let $q=\frac{n}{mn-\alpha}$, $\vec{u}^{mq}=(u_1^{mq},...,u_m^{mq})$ and $\nu=\prod_{i=1}^m{u_i^q}$. Suppose that $\vec{u}^{mq}\in A_{(1,...,1)}$ and $\nu v^{q}\in A_{\infty}$, or $u_1^{mq},...,u_m^{mq} \in A_1$ and $v^{mq} \in A_{\infty}$, then there exists a constant $C$ such that
$$\Bigg\| \frac{\mathcal M_{\alpha}(\vec f\,)(x)}{v}\Bigg\|_{L^{q, \infty}(\nu v^q)} \leq C \ \prod_{i=1}^m{\|f_i\|_{L^1(u_i)}}.$$
\end{theorem}

Note that if $\alpha=0$ then $q=\frac{1}{m}$ and we obtain Theorem \ref{muczo} for the multi(sub)linear maximal operator $\mathcal{M}$. 

\begin{Remark}
If in Theorem \ref{Max} we take $m=1$ we get that $\frac{1}{q}=1-\frac{\alpha}{n}$ and the hypothesis on the weights reduce to $u^q\in A_1$ and $v\in A_{\infty}$. Then we recover Theorem \ref{teo1.7} in the case $p=1$ for a more general class of weights $v$. The weight $u^q$ in Theorem \ref{Max} plays the role of the weight $u$ in Theorem \ref{teo1.7}.
\end{Remark}

By extrapolation arguments, we can extend this result to multilinear fractional integrals. The theorem below was essentially obtained in \cite{OP}, however, for the sake of completeness we will give a complete proof in Appendix A.

\begin{theorem}[\cite{OP}]\label{extrapolation}
Let $0<\alpha <mn$. Let $q=\frac{n}{mn-\alpha}$, $\vec{u}^{mq}=(u_1^{mq},...,u_m^{mq})\in A_{(1,...,1)}$, $v^q\in A_{\infty}$ and denote $\nu=\prod_{i=1}^m{u_i^q}$. Then there exists a constant $C$ such that
$$\Bigg\| \frac{\mathcal I_{\alpha}(\vec f\,)(x)}{v}\Bigg\|_{L^{q, \infty}(\nu v^q)} \leq C \ \Bigg\| \frac{\mathcal M_{\alpha}(\vec f\,)(x)}{v}\Bigg\|_{L^{q, \infty}(\nu v^q)}$$
\end{theorem}

Finally as a consequence of Theorem \ref{Max} and Theorem \ref{extrapolation} we obtain the main result of this paper.

\begin{theorem}\label{IMax}

Let $0<\alpha <mn$. Let $q=\frac{n}{mn-\alpha}$, $\vec{u}^{mq}=(u_1^{mq},...,u_m^{mq})$ and $\nu=\prod_{i=1}^m{u_i^q}$. Suppose that $\vec{u}^{mq}\in A_{(1,...,1)}$ and $\nu v^{q}\in A_{\infty}$, or $u_1^{mq},...,u_m^{mq} \in A_1$ and $v^{mq} \in A_{\infty}$, then there exists a constant $C$ such that
$$\Bigg\| \frac{\mathcal I_{\alpha}(\vec f\,)(x)}{v}\Bigg\|_{L^{q, \infty}(\nu v^q)} \leq C \ \prod_{i=1}^m{\|f_i\|_{L^1(u_i)}}.$$
\end{theorem}

\bigskip

The rest of the article is organized as follows. In Section 2 we recall the definition of the $A_p$ and $A_{\vec{P}}$ classes of weights. Section 3 is devoted to the proof of Theorem \ref{Max}. In Section 4, as an application of Theorem \ref{IMax}, we obtain a vector-valued extension of the mixed weighted inequalities for multilinear fractional integrals. We end this paper with an appendix, in which we give a proof of Theorem \ref{extrapolation}.

\bigskip

\textbf{Acknowledgements:} This paper constitutes a part of my doctoral thesis, under the supervision of S. Ombrosi. I would like to thank him for the support and the guidance provided during the elaboration of this paper. I would also like to thank K. Li and I. Rivera-R\'ios for reading the draft and sharing with me very useful comments and suggestions.

\section{Preliminaries}
By a weight we mean a non-negative locally integrable function defined on $\mathbb{R}^n$, such that $0<w(x)<\infty $ almost everywhere. We recall that a weight $w$ belongs to the class $A_p$, introduced by B. Muckenhoupt \cite{Mu}, $1<p<\infty$, if
$$\sup_{Q}{\Bigg( \frac{1}{|Q|} \int_Q{w(y) \ dy}\Bigg) \Bigg( \frac{1}{|Q|} \int_Q{w(y)^{1-p'} \ dy} \Bigg) ^{p-1}} <\infty, $$ where $p'$ is the conjugate exponent of $p$ defined by the equation $\frac{1}{p}+\frac{1}{p'}=1$. A weight $w$ belongs to the $A_1$ class, if there exists a constant $C$ such that $$\frac{1}{|Q|}\int_Q{w(y) \ dy}\leq C \ \inf_Q{w}. $$ Since the $A_p$ classes are increasing with respect to $p$, it is natural to define the $A_{\infty}$ class of weights by $A_{\infty}=\cup_{p\geq 1}{A_p}.$ 

\bigskip
In 2009, Lerner \textit{et al.} showed in \cite{LOPTT} that there is a way to define an analogue of the Muckenhoupt $A_p$ classes for multiple weights.
\begin{definition}
Let $m$ be a positive integer. Let $1\leq p_1,...,p_m< \infty$. We denote $p$ the number given by $\frac{1}{p}=\frac{1}{p_1}+...+\frac{1}{p_m}$, and $\vec{P}$ the vector $\vec{P}=(p_1,...,p_m)$.   
\end{definition}

\begin{definition}
Let $1\leq p_1,...,p_m<\infty $. Given $\vec{w}=(w_1,...,w_m)$, set 
$$\nu_{\vec{w}}=\prod_{i=1}^m{w_i^{\frac{p}{p_i}}}.$$
We say that $\vec{w}$ satisfies the $A_{\vec{P}}$ condition if
$$\sup_Q {\Bigg( \frac{1}{|Q|} \int_Q{\nu_{\vec{w}}}\Bigg)}^{\frac{1}{p}} \prod_{i=1}^m{\Bigg( \frac{1}{|Q|} \int_Q{w_i^{1-p'_i}}\Bigg)^{\frac{1}{p'_i}}} < \infty.$$ 
When $p_i=1$, $\Big( \frac{1}{|Q|} \int_Q{w_i^{1-p'_i}}\Big)^{\frac{1}{p'_i}}$ is understood as $(\inf_Q w_i)^{-1}$. Then we will say that $\vec{w}\in A_{(1,...,1)}$ if 
$$\sup_Q {\Bigg( \frac{1}{|Q|} \int_Q{\nu_{\vec{w}}}\Bigg)}^{\frac{1}{p}} \prod_{i=1}^m{(\inf_Q w_i)^{-1}} < \infty.$$  
\end{definition} 

The multilinear $A_{\vec{P}}$ condition has the following characterization in terms of the linear $A_p$ classes.

\begin{theorem}[\cite{LOPTT}, Theorem 3.6]\label{apmul}
Let $\vec{w}=(w_1,...,w_m)$ and $1\leq p_1,...,p_m<\infty $. Then $\vec{w}\in A_{\vec{P}}$ if and only if 

$$\left\{
w_i^{1-p_i'} \in A_{mp_i'},\ \  i=1,...,m \atop
\nu_{\vec{w}} \in A_{mp}
\right.$$ where the condition $w_i^{1-p_i'}\in A_{mp_i'}$ in the case $p_i=1$ is understood as $w_i^{\frac{1}{m}}\in A_1$.
\end{theorem}

A more general result can be found in [Lemma 3.2, \cite{LMO}].

Observe that in the particular case that every $p_i=1$ then $p=\frac{1}{m}$. By Theorem \ref{apmul}, given $\vec{w}=(w_1,...,w_m)$, we have that the following statements hold.
\begin{itemize}
\item If $\vec{w}=(w_1,...,w_m)\in A_{(1,...,1)}$ then $\nu_{\vec{w}}=w_1^{\frac{1}{m}}...w_m^{\frac{1}{m}}\in A_1$.
\item If $\vec{w}\in A_{(1,...,1)}$ then $w_i^{\frac{1}{m}}\in A_1$ for all $i=1,...,m$.
\end{itemize}

Observe that $\vec{w}\in A_{(1,...,1)}$ does not imply that $w_i\in A_1$ for every $i=1,...,m$. We can see this with a simple counterexample. Let $m=2$ and consider the weights $w_1=1$ and $w_2=\frac{1}{|x|}$. Then $\vec{w}\in A_{(1,...,1)}$, $w_1^{\frac{1}{2}}, w_2^{\frac{1}{2}}\in A_1$, $w_1\in A_1$, but $w_2\not\in A_1$.

\bigskip
 
\section{Proof of Theorem \ref{Max}}

In order to prove Theorem \ref{Max} we need the following pointwise estimate for $\mathcal M_{\alpha}$ in terms of the multilinear maximal operator $\mathcal M$. This is a multilinear version of Lemma 4 in \cite{BCP}, and to prove it, we follow a similar approach to the one that is used there.

\begin{lemma} \label{lema4}
Let $q=\frac{n}{mn-\alpha}$. Then
$$\mathcal M_{\alpha}(f_1,...,f_m)(x) \leq \mathcal M\big(f_1u_1^{1-mq},...,f_mu_m^{1-mq}\big)^{\frac{1}{mq}}(x)\prod_{i=1}^m{\Big(\int_{\mathbb{R}^n}{f_iu_i}\Big)^{\frac{\alpha}{mn}}}$$
\end{lemma}

\bigskip

\begin{proof} 
\
Let us fix $x\in \mathbb{R}^n$ and let $Q$ be a cube containing $x$. Applying H\"older's inequality with $\frac{1}{1-\frac{\alpha}{mn}}$ and $\frac{mn}{\alpha}$ we obtain

\begin{equation*}
\begin{split}
\prod_{i=1}^m{\Big( \frac{1}{|Q|^{1-\frac{\alpha}{mn}}} \int_Q {f_i} \Big)}
& = \prod_{i=1}^m{\Bigg( \frac{1}{|Q|^{1-\frac{\alpha}{mn}}} \int_Q {f_i^{1-\frac{\alpha}{mn}}f_i^{\frac{\alpha}{m}}u_i^{\frac{1}{mq}-1}u_i^{\frac{mq-1}{mq}}} \Bigg)} \\
& \leq \prod_{i=1}^m {\Bigg[ \Bigg( \frac{1}{|Q|} \int_Q{f_i u_i^{1-mq}}\Bigg)^{\frac{1}{mq}} \Bigg( \int_Q{f_iu_i}\Bigg)^{\frac{\alpha}{mn}} \Bigg]} \\
& = \prod_{i=1}^m{\Big( \frac{1}{|Q|} \int_Q {f_iu_i^{1-mq}} \Big)^{\frac{1}{mq}}} \prod_{i=1}^m{\Big( \int_{\mathbb{R}^n} {f_iu_i} \Big)^{\frac{\alpha}{mn}}}\\
& \leq \mathcal{M}\big(f_1u_1^{1-mq},...,f_mu_m^{1-mq}\big)^{\frac{1}{mq}}\prod_{i=1}^m{\Big( \int_{\mathbb{R}^n} {f_iu_i} \Big)^{\frac{\alpha}{mn}}}.
\end{split}
\end{equation*}
\\

\end{proof}

Now we have all the tools that we need to prove Theorem \ref{Max}.

\begin{proof} \textit{of Theorem \ref{Max}}
\
By applying Lemma \ref{lema4} and Theorem \ref{muczo}, we get
\begin{equation*}
\begin{split}
\lefteqn{ \nu v^q \Big\{ x\in \mathbb{R}^n : \frac{\mathcal{M}_{\alpha}(\vec f \, )(x)}{v(x)} > \lambda \Big\}^{\frac{1}{q}}} \\         
& \, \, \leq u_1^q...u_m^q v^q \Big\{ x\in \mathbb{R}^n : \frac{\mathcal{M}(f_1 u_1^{1-mq},...,f_m u_m^{1-mq})^{\frac{1}{mq}}(x)}{v(x)} > \frac{\lambda}{\prod_{i=1}^m{(\int_{\mathbb{R}^n}{f_iu_i})^{\frac{\alpha}{mn}}}} \Big\}^{\frac{1}{q}} \\
& = u_1^q...u_m^q v^q \Big\{ x\in \mathbb{R}^n : \frac{\mathcal{M}(f_1 u_1^{1-mq},...,f_m u_m^{1-mq})(x)}{v^{mq}(x)} > \Big( \frac{\lambda}{\prod_{i=1}^m{(\int_{\mathbb{R}^n}{f_iu_i})^{\frac{\alpha}{mn}}}} \Big)^{mq} \Big\}^{\frac{1}{q}} \\
& = (u_1^{mq})^{\frac{1}{m}}...(u_m^{mq})^{\frac{1}{m}} (v^{mq})^{\frac{1}{m}} \Big\{ x\in \mathbb{R}^n : \frac{\mathcal{M}(f_1 u_1^{1-mq},...,f_m u_m^{1-mq})(x)}{v^{mq}(x)} > \Big( \frac{\lambda}{\prod_{i=1}^m{(\int_{\mathbb{R}^n}{f_iu_i})^{\frac{\alpha}{mn}}}} \Big)^{mq} \Big\}^{m \frac{1}{mq}} \\
& \leq \frac{C}{\lambda} \prod_{i=1}^m{\Big(\int_{\mathbb{R}^n}{f_iu_i}\Big)^{\frac{\alpha}{mn}}} \prod_{i=1}^m{\Big(\int_{\mathbb{R}^n}{f_iu_i^{1-mq}u_i^{mq}} \Big)^{\frac{1}{mq}}}\\
& = \frac{C}{\lambda} \prod_{i=1}^m{\Big(\int_{\mathbb{R}^n}{f_iu_i} \Big)} \\
& = \frac{C}{\lambda} \prod_{i=1}^m{\|f_i\|_{L^1{u_i}}}.
\end{split}
\end{equation*}
\\

\end{proof}

\section{A vector-valued extension of Theorem \ref{IMax}}

Recently in \cite{CMO} D. Carando, M. Mazzitelli and S. Ombrosi obtained a generalization of the Marcinkiewicz-Zygmund inequalities to the context of multilinear operators. We recall one of the results in \cite{CMO} that extends previously known results from \cite{GM} and \cite{BPV}. 

\begin{theorem}[\cite{CMO}]\label{CMO}
Let $0<p, q_1, \dots, q_m < r<2$ or $r=2$ and $0<p, q_1, \dots, q_m < \infty$ and, for each $1 \leq i \leq m$, consider $\{f^i_{k_i}\}_{k_i} \subset L^{q_i}(\mu_i)$.  And Let $S$ be a multilinear operator such that  
 $S\colon L^{q_1}(\mu_1) \times \cdots \times L^{q_m}(\mu_m) \to L^{p,\infty}(\nu)$, then, there exists a constant $C>0$ such that 
 
 \begin{equation}\label{weak MZ multilineal p,q>0}
\left\| \left( \sum_{k_1, \dots, k_m} |S(f^1_{k_1}, \dots, f^m_{k_m})|^r \right)^{\frac{1}{r}}  \right\|_{L^{p, \infty}(\nu)} \leq C  \|S\|_{weak} \prod_{i=1}^m \left\| \left( \sum_{k_i} |f^i_{k_i}|^{r} \right)^{\frac{1}{r}} \right\|_{L^{q_i}(\mu_i)}.
\end{equation}

\end{theorem}

As a consequence of this theorem and Theorem \ref{IMax}  we obtain the following mixed weighted vector valued inequality for a multilinear fractional operator $\mathcal{I}_{\alpha}$. 

\begin{cor}\label{ext vec} 
Let $S(\vec{f \,})=\frac{\mathcal{I}_{\alpha}(\vec{f\,})}{v}$, where $\mathcal{I}_{\alpha}$ is a multilinear fractional operator. Let $q=\frac{n}{mn-\alpha}$, $\vec{u}^{mq}=(u_1^{mq},...,u_m^{mq})$ and $\nu=\prod_{i=1}^m{u_i^q}$. Suppose that $\vec{u}^{mq}\in A_{(1,...,1)}$ and $\nu v^{q}\in A_{\infty}$, or $u_1^{mq},...,u_m^{mq} \in A_1$ and $v^{mq} \in A_{\infty}$.
For each $1 \leq i \leq m$, consider $\{f^i_{k_i}\}_{k_i} \subset L^{1}(u_i)$. Then, there exists a constant $C>0$ such that 

\begin{equation}\label{ext vect}
\left\| \left( \sum_{k_1, \dots, k_m} |S(f^1_{k_1}, \dots, f^m_{k_m})|^r \right)^{\frac{1}{r}}  \right\|_{L^{q, \infty}(\nu v^q)} \leq C \prod_{i=1}^m \left\| \left( \sum_{k_i} |f^i_{k_i}|^{r} \right)^{\frac{1}{r}} \right\|_{L^{1}(u_i)}.
\end{equation}
\end{cor}

Observe that under the hypothesis of Corollary \ref{ext vec}, $S$ satisfies $S\colon L^{1}(u_1) \times \cdots \times L^{1}(u_m) \to L^{q,\infty}(\nu v^q)$. So we are under the hypothesis of Theorem \ref{CMO}.

\bigskip

\section{Appendix A. Proof of Theorem \ref{extrapolation}}
In order to settle Theorem \ref{extrapolation} we will need two known results. The first one is due to K. Moen \cite{M}.
\begin{theorem}[\cite{M}, Theorem 3.1]\label{teomoen}
Suppose that $0<\alpha<mn$, then for every $w\in A_{\infty}$ and all $0<s<\infty$ we have
$$\int_{\mathbb{R}^n}{|\mathcal{I}_{\alpha}\vec{f}(x)|^s \ w(x)\ dx}\leq C\int_{\mathbb{R}^n}{\mathcal{M}_{\alpha}\vec{f}(x)^s\ w(x)\ dx}$$ for all functions $\vec{f}$ with $f_i$ bounded with compact support.
\end{theorem}

\bigskip
The second result we will rely upon is due to D. Cruz-Uribe, J.M. Martell and C. P\'erez \cite{CMP}.
\begin{theorem}[\cite{CMP}, Theorem 1.7]\label{teocmp}
Given a family $\mathcal{F}$ of pairs of functions that satisfies that there exists a number $p_0$, $0<p_0<\infty$ such that for all $w\in A_{\infty}$
$$\int_{\mathbb{R}^n}{f(x)^{p_0}\ w(x) \ dx} \leq C \int_{\mathbb{R}^n}{g(x)^{p_0}\ w(x) \ dx},$$ for all $(f,g)\in \mathcal{F}$ such that the left hand side is finite, and with $C$ depending only on $[w]_{A_{\infty}}$. Then, for all weights $u$, $v$ such that $u\in A_1$ and $v\in A_{\infty}$ we have that
$$\|fv^{-1}\|_{L^{1,\infty}(uv)}\leq C \|gv^{-1}\|_{L^{1,\infty}(uv)} \  \, \, (f,g)\in \mathcal{F}.$$
\end{theorem}

\bigskip

Having those results at our disposal we proceed as follows.
First of all observe that if $\vec{u}^{mq}=(u_1^{mq},...,u_m^{mq}) \in A_{(1,...,1)}$ then $\nu = u_1^q...u_m^q \in A_1$. 

Then, by Theorem \ref{teomoen} and Theorem \ref{teocmp},
\begin{equation*}
\begin{split}
\Bigg\| \frac{\mathcal{I}_{\alpha}(\vec{f\,})}{v}\Bigg\|_{L^{q, \infty}(\nu v^q)}^q
& = \sup_{\lambda>0} \ \lambda^q \Big(\nu v^q \ \{ x\in {\mathbb{R}}^n:\Bigg|\frac{\mathcal{I}_{\alpha}(\vec{f\,})(x)}{v(x)}\Bigg|>\lambda \}\Big)\\
& = \sup_{\lambda>0} \ \lambda^q \Big(\nu v^q \ \{ x\in {\mathbb{R}}^n:\Bigg|\frac{T(\vec{f\,})(x)}{v(x)}\Bigg|^q>\lambda^q \}\Big) \\
& = \sup_{t>0} \ t \ \Big(\nu v^q \ \{ x\in {\mathbb{R}}^n:\Bigg|\frac{T(\vec{f\,)}(x)}{v(x)}\Bigg|^q>t\}\Big) \\
& = \Bigg\| \Bigg(\frac{\mathcal{I}_{\alpha}(\vec{f\,})}{v}\Bigg)^q \Bigg\|_{L^{1, \infty}(\nu v^q)}\\
& \leq C \Bigg\| \Bigg(\frac{\mathcal{M}_{\alpha}(\vec{f\,})}{v}\Bigg)^q \Bigg\|_{L^{1, \infty}(\nu v^q)}\\
& = \sup_{\lambda>0} \ \lambda \Big(\nu v^q \ \{ x\in {\mathbb{R}}^n:\Bigg( \frac{\mathcal{M}_{\alpha}(\vec{f\,})(x)}{v(x)}\Bigg)^q>\lambda \}\Big)\\
& = \sup_{t>0} \ t^q \Big(\nu v^q \ \{ x\in {\mathbb{R}}^n:\Bigg( \frac{\mathcal{M}_{\alpha}(\vec{f\,})(x)}{v(x)}\Bigg) ^q>t^q \}\Big) \\
& = \sup_{t>0} \ t^q \ \Big(\nu v^q \ \{ x\in {\mathbb{R}}^n: \frac{\mathcal{M}_{\alpha}(\vec{f\,)}(x)}{v(x)}>t\}\Big) \\
& = \Bigg\| \frac{\mathcal{M}_{\alpha}(\vec{f\,})}{v} \Bigg\|_{L^{q, \infty}(\nu v^q)}^q.
\end{split}
\end{equation*}

\end{document}